\documentclass[conference]{IEEEtran}
\IEEEoverridecommandlockouts
\usepackage{graphicx}
\usepackage{graphics}
\usepackage{subfig}
\usepackage{amsmath}
\usepackage{algorithm}
\usepackage{algorithmic}
\usepackage{cite}
\usepackage{amsthm}
\usepackage{amssymb}
\usepackage[table]{xcolor}
\usepackage{multicol}
\usepackage{bm}
\usepackage{tabularx}
\usepackage{booktabs}
\usepackage{makecell}
\usepackage{multirow}
\usepackage{threeparttable}
\usepackage{blindtext}
\usepackage{scrextend}
\addtokomafont{labelinglabel}{\sffamily}

\allowdisplaybreaks[4]

\theoremstyle{definition}

\theoremstyle{remark}

\theoremstyle{proposition}

\usepackage[font=small,labelsep=period]{caption}
\definecolor{ColorName}{rgb}{0 0 1}
\definecolor{Red}{rgb}{1 0 0}

\makeatletter
\DeclareRobustCommand{\bxz}{\ifmmode \mathbxz
	\else
	\leavevmode\unskip\penalty9999 \hbox{}\nobreak\hfill
	\quad\hbox{\bxzsymbol}%
	\fi
}
\newcommand{\mathbxz}{\quad\hbox{\bxzsymbol}}

\providecommand{\bxzsymbol}{\fbox{\footnotesize}}
\providecommand{\fooname}{Proof of Theorem \ref{Apro. constant}}
\makeatother
\begin{document}

\title{Energy-Efficient Data Offloading for \\Earth Observation Satellite Networks}

\author{

\IEEEauthorblockN{
	Lijun He\IEEEauthorrefmark{1},
	Ziye Jia\IEEEauthorrefmark{2},
	Juncheng Wang\IEEEauthorrefmark{3},
	Feng Wang\IEEEauthorrefmark{4},
	Erick Lansard\IEEEauthorrefmark{5},
	and Chau Yuen \IEEEauthorrefmark{5}\\
	}
	\IEEEauthorblockA{\IEEEauthorrefmark{1} School of Software, Northwestern Polytechnical University, China\\
		\IEEEauthorrefmark{2}
Ministry of Industry and Information Technology, Nanjing University of Aeronautics and Astronautics, China\\
\IEEEauthorrefmark{3}
Department of Computer Science, Hong Kong Baptist University, China
\\
\IEEEauthorrefmark{4}
Information Systems Technology and Design, Singapore University of Technology and Design, Singapore\\
\IEEEauthorrefmark{5}
School of Electrical and Electronics Engineering, Nanyang Technological University, Singapore
\\
Email: lijunhe@nwpu.edu.cn, jiaziye@nuaa.edu.cn, jcwang@comp.hkbu.edu.hk, feng2\_wang@sutd.edu.sg,\\
\{chau.yuen, erick.lansard\}@ntu.edu.sg}

\thanks{This work was supported in part by the National Natural Science Foundation of China under Grant 62201463, in part by Natural Science Basic Research Program of Shaanxi Province under Grant 2022JQ-615, in part by Basic Research Programs of Taicang under Grant TC2021JC25.}

}


\maketitle
\begin{abstract}
In Earth Observation Satellite Networks (EOSNs) with a large number of battery-carrying satellites, proper power allocation and task scheduling are crucial to improving the data offloading efficiency. As such, we jointly optimize power allocation and task scheduling to achieve energy-efficient data offloading in EOSNs, aiming to balance the objectives of reducing the total energy consumption and increasing the sum weights of tasks. First, we derive the optimal power allocation solution to the joint optimization problem when the task scheduling policy is given. Second, leveraging the conflict graph model, we transform the original joint optimization problem into a maximum weight independent set problem when the power allocation strategy is given. Finally, we utilize the genetic framework to combine the above special solutions as a two-layer solution for the joint optimization problem. Simulation results demonstrate that our proposed solution can properly balance the sum weights of tasks and the total energy consumption, achieving superior system performance over the current best alternatives.
\end{abstract}


\IEEEpeerreviewmaketitle

\section{Introduction}
Driven by the fast proliferation of remote sensing applications, such as environment monitoring, meteorology, and natural disaster surveillance, Earth Observation Satellite Networks (EOSNs) experience an explosive growth in remote sensing data \cite{Du2016Cooperativ,He2022JOTS,Zhou2023AINI}. There arises an urgent need to offload the collected remote sensing data from EOSNs to Ground Stations (GSs) for further data processing. Different from traditional ground wireless networks, EOSNs belong to Ad Hoc networks with dynamic network topologies \cite{Jiandong2002CEMM}. This means that data offloading in EOSNs occurs only within limited Transmission Time Windows (TTWs) between Earth Observation Satellites (EOSs) and GSs due to the high-speed motion of EOSs. The key to efficient data offloading in EOSNs is how to schedule data offloading tasks within limited TTWs to maximize the desired network performance\cite{Deng2018TPTS}.

From the mathematical perspective, the data offloading problem in EOSNs is similar to the Parallel Machine Scheduling Problem with Time Windows (PMSPTW)\cite{Li2023TS}. Recently, a majority of works studied the PMSPTW to select a set of optimal GSs for Low Earth Orbit (LEO) satellites to schedule Satellite Ground Links (SGLs) while satisfying the time window constraints \cite{Zufferey2015treeSearch,Liu2019DP,Chen2012PSO,ZHANG2018ACO,ZHANG2022GA,Qu2023DRL,Liang2023DRL}. The scheduling algorithms proposed in these works are either deterministic \cite{Zufferey2015treeSearch,Liu2019DP} or non-deterministic\cite{Chen2012PSO,ZHANG2018ACO,ZHANG2022GA,Qu2023DRL,Liang2023DRL}. The deterministic algorithms use mathematical tools, such as tree search\cite{Zufferey2015treeSearch} and dynamic programming\cite{Liu2019DP}, to find an optimal solution under mild assumptions. However, these deterministic algorithms are time-consuming especially for large-scale problems and thus are not suitable for practical EOSNs. To accommodate large-scale EOSNs, some non-deterministic algorithms have been proposed to find approximation solutions with lower computation costs. They are mainly based on meta-heuristic methods, such as particle swarm optimization\cite{Chen2012PSO}, ant colony optimization\cite{ZHANG2018ACO}, and genetic algorithm \cite{ZHANG2022GA}. Recently, some attempts have been made to utilize Deep Reinforcement Learning (DRL) for solving satellite scheduling problems. For example, the PMSPTW was transformed into an Assignment Problem (AP) and a Single Antenna Scheduling Problem (SASP) \cite{Qu2023DRL}. Then, the DRL method and a heuristic scheduling method were proposed to solve the AP and the SASP, respectively. In \cite{Liang2023DRL}, a DRL-based algorithm was proposed for a rapid satellite range rescheduling problem to process real-time emergency events.

However, the above works ignored the energy management of battery-carrying LEO satellites powered by the solar panels, resulting in data offloading interruptions due to the overuse of energy. As such, energy consumption control is crucial for efficient data offloading in EOSNs\cite{zhou2020mission}.
To improve data offloading efficiency, some researchers studied energy-efficient offloading for various satellite networks, such as land mobile satellite systems\cite{An2018Power}, cognitive satellite terrestrial networks\cite{Ruan2019EE}, small satellite cluster networks \cite{Zhang2022jointddrm}, and satellite-terrestrial networks \cite{Gao2023ECOS}. In these works, each offloading task is represented as a set of equal packets, and then energy-efficient scheduling strategies are proposed to download these packets. In this case, the decision variables for scheduling packets form a huge solution space, increasing the difficulty of solving the problem. Therefore, these works \cite{An2018Power,Ruan2019EE,Gao2023ECOS,Zhang2022jointddrm} are unsuitable for EOSNs to offload the tasks with a large amount of data.


\begin{figure}[!t]
	\centering	\includegraphics[scale=0.55]{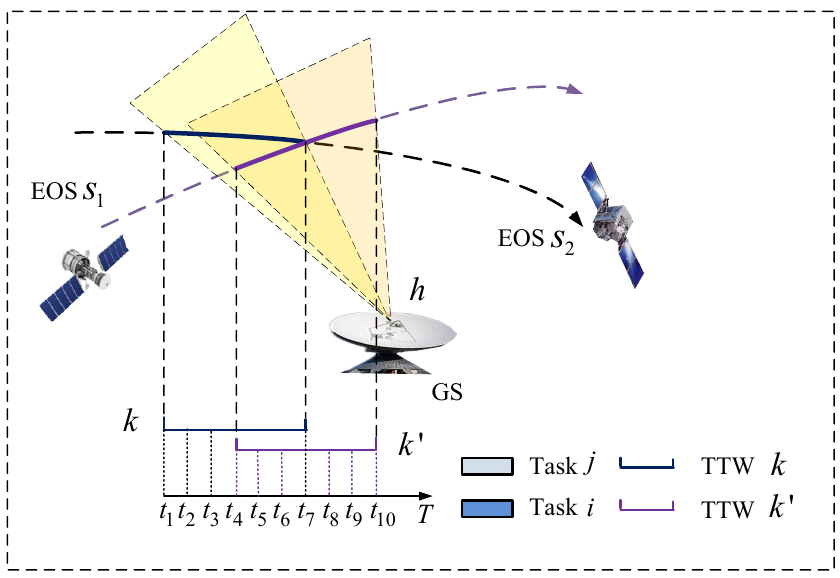}
	\caption{Illustration of data offloading in EOSNs.}
	\label{Fig:SystemScenario}
\end{figure}

In contrast to previous works \cite{An2018Power,Ruan2019EE,Gao2023ECOS,Zhang2022jointddrm}, in this work, we character the main features of data offloading in real EOSNs from the following three aspects. First, we consider a more realistic data offloading scenario with variable-size tasks. Second, our proposed joint optimization framework incorporates the time window constraints of data offloading in EOSNs, which captures the intermittent characteristic of SGLs. Third, we develop an efficient yet fast two-layer optimization solution, termed Energy-efficient Data Offloading (EDO), to minimize the total energy consumption while maximizing the sum weights of tasks.

\section{System Model and Problem Formulation} \label{sec:SystemModel}

We consider a downlink EOSN consisting of a set $\mathcal{S}\!=\!\{1,...,S\}$ of EOSs and a set of GSs. We discretize the whole scheduling horizon into $T$ time slots to obtain a set of time slots, denoted by $\mathcal{T}\!=\!\{1,...,T\}$. We denote by time slot $t$ the time interval of $[t,t+1]$. We assume that each EOS and each GS is respectively equipped with one transmit antenna and one receiver antenna. Let $\mathcal{H}=\{1,...,H\}$ denote the set of data receiver antennas. As shown in Fig. \ref{Fig:SystemScenario}, EOSs cycle their individual orbits at high speed, thereby yielding a set of short contact TTWs. Let $\mathcal{K}_{s,h}$ denote the set of TTWs between each EOS $s\in\mathcal{S}$ and each antenna $h\in\mathcal{H}$. Let $\mathcal{K}_{s,\mathcal{H}}=\{\mathcal{K}_{s,h}, h\in \mathcal{H}\}$. We use 4-tuples $(s,h,a_k,b_k)$ to represent each TTW $k\in\mathcal{K}_{s,h}$, where $a_k$ and $b_k$  are the beginning and end times, respectively. We denote by $\mathcal{J}=\{1,...,J\}$ the set of data offloading tasks. Let $s(j)$ index the EOS that stores the data of task $j$.


\subsection{Channel Model}
We denote by $R_{jk}^{\text{SGL}}$ the achievable data rate of SGLs within TTW $k$ in bits/s for any task $j$, given by
\begin{align}\label{Def:RtoSGLs}
R_{jk}^{\text{SGL}}=B_c\log_2(1+\text{SNR}_{jk}),~~\forall j\in \mathcal{J},k\in\mathcal{K}_{s(j),\mathcal{H}},
\end{align}
where $B_c$ is the bandwidth and $\text{SNR}_{jk}$ is the signal-to-noise ratio of offloading task $j$ within TTW $k$. We evaluate the $\text{SNR}_{jk}$ of SGLs \cite{Golkar2015FSSP} as follows:
\begin{align}\label{Def:SNR}
\text{SNR}_{jk} = \frac{P_{s(j)}G^{\text{tran}}_{s(j)}G^{\text{rec}}_hL_fL^k_l}{N},
\end{align}
where $G^{\text{tran}}_{s(j)}$ represents the transmit antenna gain of EOS $s(j)$, $G^{\text{rec}}_h$ denotes the gain of receiver antenna $h$, $L_f$ is the free space loss, $N$ is the noise power, and $L_l^{k}$ is the total path loss for TTW $k$. The total path loss in dB is calculated by $L_l^{k}=L^k_b+L^k_g+L^k_s+L^k_e$, where $L^k_b$ denotes the basic path loss in dB, $L^k_g$ represents the attenuation due to atmospheric gasses in dB, $L^k_s$ indicates the attenuation due to either ionospheric or tropospheric scintillation in dB, and $L^k_e$ is the building entry loss in dB for TTW $k$ \cite{3gppDocu}. In particular, $P_{s(j)}$ represents the transmit power of EOS $s(j)$ allocated to task $j$, which is subject to the following transmit power range of each EOS:
\begin{align}
\text{C1}:~0\le P_{s(j)} \le P_{s(j)}^{\text{max}},~~\forall j\in\mathcal{J}.\label{P0:C5}
\end{align}

Furthermore, we let $R_{k}^{\text{req}}$ represent the rate requirement of TTW $k$ and introduce the following constraint to indicate the rate requirement of each TTW:
\begin{align}
\text{C2}:~R_{jk}^{\text{SGL}} \ge R_{k}^{\text{req}}, ~~\forall j\in \mathcal{J},k\in\mathcal{K}_{s(j),\mathcal{H}}.\label{P0:C4}
\end{align}

\begin{figure}[!t]
	\centering	\includegraphics[scale=0.83]{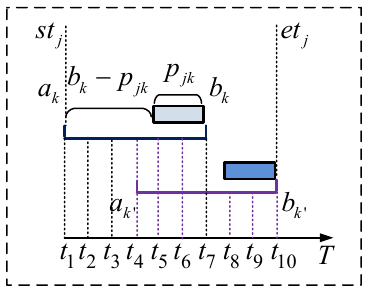}
	\caption{Illustration of the set of feasible time slots within TTW $k$.}
	\label{Fig:Timeslot}
\end{figure}

\subsection{Task Model}
Let $D_j$ and $w_j$ denote the data amount of task $j$ and the weighted value of task $j$, respectively. We denote by $p_{jk}$ the processing time, i.e., the transmitting time, of task $j$ within TTW $k$, which can be computed by
\begin{align}\label{Def:pjk}
	p_{jk}=\frac{D_j}{R_{jk}^{\text{SGL}}}.
\end{align}
Let $st_j$ and $et_j$ be the earliest time and the latest time to start offloading task $j$, respectively.
We use 6-tuples $(s(j),D_j,w_j,p_{jk},st_j,et_j)$ to represent each task $j\in\mathcal{J}$. Furthermore, we use a decision variable $x^t_{jk}\in\{0,1\}$ to represent the offloading strategy, such that $x^t_{jk}=1$ indicates that the data of task $j$ is offloaded in time slot $t$ within TTW $k$; and $x^t_{jk}=0$ otherwise. Thus, we introduce the following binary constraint
\begin{align}
\text{C3}:~x^t_{jk}\in\{0,1\}, ~~\forall j\in \mathcal{J},k\in\mathcal{K}_{s(j),\mathcal{H}},t\in\mathcal{F}(j),\label{P0:C6}
\end{align}
to indicate whether task $j$ starts to offload in time slot $t$ within TTW $k$, or not. Wherein, $\mathcal{F}(j)$ represents the set of all feasible time slots for task $j$, given by
\begin{align}
\mathcal{F}(j)=\bigcup\limits_{k\in\mathcal{K}_{s(j),\mathcal{H}}}\mathcal{F}(j,k),
\end{align}
with $\mathcal{F}(j,k)=[a_k,b_k-p_{jk}]\cap[st_j,et_j]$ indicating the set of feasible slots within TTW $k$ to start offloading the data of task $j$. As an example for illustration in Fig. \ref{Fig:Timeslot}, we have $\mathcal{F}(j,k)=\{t_1,t_2,...,t_5\}$.

We further introduce the following constraint to guarantee that each task $j$ is offloaded at most once:
\begin{align}
\text{C4}:~\sum_{t\in\mathcal{F}(j)}x^{t}_{jk}\le 1,~~\forall j\in\mathcal{J}.\label{P0:C1}
\end{align}
In addition, we denote by $E_j$ the energy consumption for offloading task $j$, which is calculated by
\begin{align}\label{Def:EJ}
E_j = \sum_{t\in\mathcal{F}(j)}x^{t}_{jk}P_{s(j)}p_{jk}.
\end{align}

\subsection{Time Window Model}
We define the set of occupying time slots for antenna $h$ to receive the data of task $j$ in time slot $t$ as
\begin{align*}
&\Phi(j,h,t) = \bigcup\limits_{k\in\mathcal{K}_{s(j),h}}\left[t,t+p_{jk}-1\right]\cap[a_k,b_k]\cap[st_j,et_j].
\end{align*}
We can see from the example in Fig. \ref{Fig:Timeslot} that $\Phi(j,h,t_5)=\{t_5,t_6\}$. As such, the following constraint reveals that any receiver antenna receives the data of at most one task in any time slot:
\begin{align}
\text{C5}:~\sum\limits_{j\in \mathcal{J}}\sum\limits_{\tau\in\Phi(j,h,t)}x^{\tau}_{jk}\le 1,~~\forall h\in\mathcal{H}, t\in \mathcal{F}(j).\label{P0:C2}
\end{align}
Similarly, we define the set of occupying feasible time slots for EOS $s$ to offload the data of task $j$ in time slot $t$ as follows:
\begin{align*}
&\Theta(j,s,t) = \bigcup\limits_{k\in\mathcal{K}_{s(j),\mathcal{H}}}\left[t,t+p_{jk}-1\right]\cap[a_k,b_k]\cap[st_j,et_j].
\end{align*}
We further use the following constraint to indicate that any EOS offloads the data of at most one task in any time slot:
\begin{align}
\text{C6}:~\sum\limits_{j\in \mathcal{J}}\sum\limits_{\tau\in\Theta(j,s,t)}x^{\tau}_{jk}\le 1,~~\forall s\in\mathcal{S},t\in \mathcal{F}(j).\label{P0:C3}
\end{align}

\subsection{Problem Formulation}


To balance the two conflicting objectives of minimizing the total energy consumption and maximizing the sum weights of tasks, we consider their normalized weighted sum and formulate the following mixed integer programming problem:
\begin{align}
	\textbf{P0}(\bm{x},\bm{P})\text{:}~~&\min_{\bm{x},\bm{P}}~~\frac{\lambda}{E_{\text{max}}^{\text{total}}}\sum_{j \in \mathcal{J}}E_j-\sum_{j \in \mathcal{J}}\sum_{t\in\mathcal{F}(j)}\frac{1-\lambda}{W_{\text{max}}}w_jx^{t}_{jk}\nonumber\\
	&\text{s.t.}~~\text{C1-C6},
\end{align}
where $\bm{x}=\{x^t_{jk}\}$, $\bm{P}=\{P_{s(j)}\}$, $W_{\text{max}}=\sum_{j \in \mathcal{J}}w_j$ denotes the sum weights of all the tasks, $E_{\text{max}}^{\text{total}}$ represents the sum of the maximization energy consumption of all tasks, and $\lambda\in[0,1)$ is a tuning weighted coefficient to obtain the desired total energy consumption and sum weights of tasks. Note that the distinct characteristic of the data offloading in EOSNs is mainly reflected in the time window constraints of C5 and C6.




\section{Proposed Solution}
In this section, we first study two cases for $\textbf{P0}(\bm{x},\bm{P})$ to obtain two special solutions and then utilize the genetic framework to combine them as an efficient two-layer solution.

\subsection{Power Allocation for Fixed Task Scheduling} \label{Subsec:PowerAllocation}
With a given task scheduling decision $\bm{x}$, $\textbf{P0}(\bm{x},\bm{P})$ degrades into the total energy consumption minimization problem:
\begin{align*}
\textbf{P1}(\bm{P}):~~\eta(\bm{x})~=~&\min_{\bm{P}}~~\frac{\lambda}{E_{\text{max}}^{\text{total}}}\sum_{j \in \mathcal{J}}E_j-\mathcal{U}\\
&\text{s.t.}~~\text{C1-C2},
\end{align*}
where $\mathcal{U}$ is a constant, given by
\begin{align}
\mathcal{U} = \sum_{j \in \mathcal{J}}\sum_{t\in\mathcal{F}(j)}\frac{1-\lambda}{W_{\text{max}}}w_jx^{t}_{jk}.
\end{align}

We further explore the unique structure of $\textbf{P1}(\bm{P})$ to obtain the optimal power allocation solution. We observe that $\textbf{P1}(\bm{P})$ can be decomposed into $J$ independent subproblems and one for each task $j$, given by
\begin{align*}
\textbf{P2}(\bm{P})\text{:}~~&\min_{P_{s(j)}}~~\frac{\lambda}{E_{\text{max}}^{\text{total}}}E_j\\
&\text{s.t.}~~\text{C1-C2}.
\end{align*}
In line with (\ref{Def:RtoSGLs}) and (\ref{Def:SNR}), we obtain
\begin{align}\label{Def:Rjk}
R_{jk}^{\text{SGL}}=B_c\log_2(1+\frac{P_{s(j)}G^{\text{tran}}_{s(j)}G^{\text{rec}}_hL_fL_l^{k}}{N}),\nonumber\\
~~\forall j\in\mathcal{J}, k\in\mathcal{K}_{s(j),\mathcal{H}}.
\end{align}
Combining (\ref{Def:pjk}) and (\ref{Def:EJ}), we have
\begin{align}\label{Def:E_jnew}
E_j = \sum_{t\in\mathcal{F}(j)}\frac{x^{t}_{jk}P_{s(j)}D_j}{R_{jk}^{\text{SGL}}},~~\forall j\in \mathcal{J}.
\end{align}
We combine (\ref{Def:Rjk}) with (\ref{Def:E_jnew}) to obtain:
\begin{align}
E_j = \frac{\alpha P_{s(j)}}{\log_2(1+\beta P_{s(j)})},~~\forall j\in \mathcal{J},
\end{align}
where $\alpha=\sum_{t\in\mathcal{F}(j)}\frac{x^{t}_{jk}D_j}{B_c}$ and $\beta=\frac{G^{\text{tran}}_{s(j)}G^{\text{rec}}_hL_fL_l^{k}}{N}$. To study the properties of $E_j$, we first introduce a new variable $\bm{y}=\{y_j\}$ with each element satisfying
\begin{align}
y_j = \frac{1}{\log_2(1+\beta P_{s(j)})},~~\forall j\in \mathcal{J}.
\end{align}
We further use $\bm{y}$ to recast {$\textbf{P2}(\bm{P})$} as the following new form:
\begin{align*}
\textbf{P3}(\bm{y})\text{:}~~&\min_{y_j}~~\frac{\alpha}{\beta}y_j(2^{\frac{1}{y_j}}-1)\\
&\text{s.t.}~~f(P^{\text{max}}_{s(j)})\le y_j\le\frac{B_c}{R_k^{\text{req}}},
\end{align*}
where $f(P^{\text{max}}_{s(j)})=\frac{1}{\log_2(1+\beta P^{\text{max}}_{s(j)})}$ and $g(y_j)=y_j(2^{\frac{1}{y_j}}-1)$. Next, we can obtain the first-order and second-order derivation of $g(y_j)$ as follows:
\begin{align}
&\triangledown g(y_j) = 2^{\frac{1}{y_j}}\left(1-\frac{\ln 2}{y_j}\right)-1,\\
&\triangledown^2 g(y_j) = 2^{\frac{1}{y_j}}\frac{(\ln 2)^2}{(y_j)^3}.
\end{align}

For any $y_j$ subject to $f(P^{\text{max}}_{s(j)})\le y_j\le\frac{B_c}{R_k^{\text{req}}}$, we have $\triangledown^2 g(y_j)>0$, which indicates that $g(y_j)$ is a convex and unimodal function.
Thus, the optimal solution to ${\textbf{P3}}(\bm{y})$ is either $f(P^{\text{max}}_{s(j)})$, $\frac{B_c}{R_k^{\text{req}}}$, or the extreme point of $g(y_j)$, depending on the value of $\hat{y}_j$, given by

\begin{align}
y^{*}_j=\left\{
\begin{array}{lcl}
f(P^{\text{max}}_{s(j)})      &     & {\hat{y}_j<f(P^{\text{max}}_{s(j)})},\\
\hat{y}_j   &     & {f(P^{\text{max}}_{s(j)})\le \hat{y}_j<\frac{B_c}{R_k^{\text{req}}}},\\
\frac{B_c}{R_k^{\text{req}}}     &     & {\frac{B_c}{R_k^{\text{req}}}\le \hat{y}_j},
\end{array} \right.
\end{align}
where $\hat{y}_j$ satisfies
$\triangledown g(\hat{y}_j) = 2^{\frac{1}{\hat{y}_j}}\left(1-\frac{\ln 2}{\hat{y}_j}\right)-1=0$. In particular, we can use the bisection method to obtain $\hat{y}_j,\forall j\in \mathcal{J}$. Therefore, we obtain the optimal solution to $\textbf{P1}(\bm{P})$:
\begin{align}\label{Def:P*_s(j)}
\bm{P}^*=\{P^*_{s(j)}|P^*_{s(j)}=\frac{1}{\beta}(2^{y^*_j}-1),~\forall j\in \mathcal{J}\}.
\end{align}
\subsection{Task Scheduling for Fixed Power Allocation} \label{Subsec:TaskScheduingFixedPower}
With a given power allocation solution $\bm{P}$, $\textbf{P0}(\bm{x},\bm{P})$ degrades into the following task scheduling problem:
\begin{align}
\textbf{P4}(\bm{x})\text{:}~~& \max_{\bm{x}}~~\sum_{j \in \mathcal{J}}\sum_{t\in\mathcal{F}(j)}\left(\frac{1-\lambda}{W_{\text{max}}}w_j-\frac{\lambda}{E_{\text{max}}^{\text{total}}}P_{s(j)}p_{jk}\right)x^{t}_{jk}\nonumber\\
&\text{s.t.}~~\text{C3-C6}.\nonumber
\end{align}

\begin{figure}[!t]
	\centering
\includegraphics[scale=0.65]{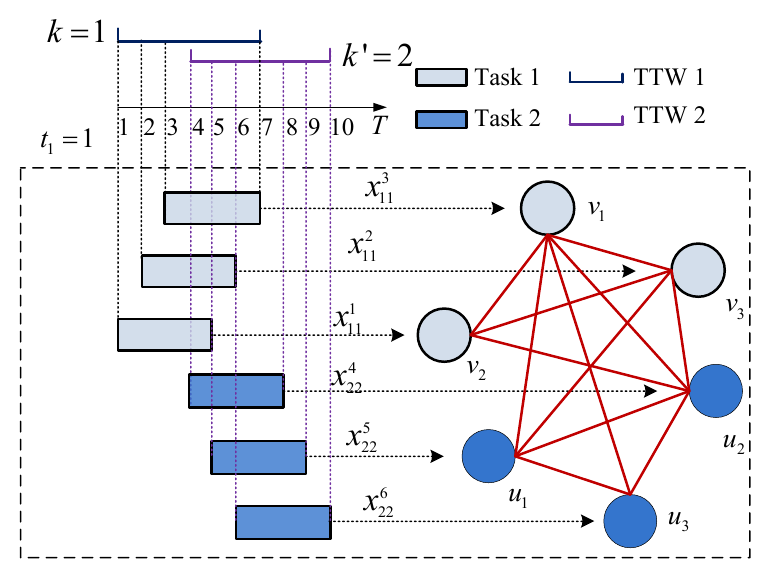}
	\caption{The construction of conflict graph.}
	\label{Fig:ConflictGraph}
\end{figure}

To solve $\textbf{P4}(\bm{x})$ efficiently, we leverage a conflict graph model $G(V,E)$ with $V$ and $E$ respectively representing the vertex set and edge set, to depict the conflict relations among variables. We take the following two steps to construct conflict graph $G(V,E)$. First, we construct vertex set $V$ by representing any variable $x^{t}_{jk}$ as a vertex $v$ one by one. As such, we establish one to one mapping between $x^{t}_{jk}$ and $v$. Second, from the constraints of $\textbf{P4}(\bm{x})$, we can check whether any two variables $x^{t}_{jk}$ (denoted by $u$) and $x^{t'}_{j'k'}$ (denoted by $v$) take the value of one at the same time. If not, we can add one edge across the vertexes $u$ and $v$. We perform this procedure on all combinations on $\bm{x}$ to constitute vertex set $E$. We give an example for illustrating the conflict graph construction in Fig.~\ref{Fig:ConflictGraph}. By constructing the conflict graph model $G(V,E)$, we can equivalently transform $\textbf{P4}(\bm{x})$ into a Maximum Weight Independent Set (MWIS) problem .

\subsection{Joint Power Allocation and Task Scheduling}\label{Subsec:JTPA}

\begin{algorithm}[!t]\caption{Energy-efficient Data Offloading (EDO)}\label{alg:EDO}
	\begin{algorithmic}[1]
		\REQUIRE $P^{\text{max}}_{s(j)}$, generation number $L$, population size $U$.\\
        \STATE {Obtain $\bm{P}^*$ by calling the solution in Section \ref{Subsec:PowerAllocation}. }
    \STATE{Utilize $\{P^{\text{max}}_{s(j)}\}$ and $\bm{P}^*$ to obtain $\mathcal{P}_{\text{Discret}}$.}
    \STATE{Use $\mathcal{P}_{\text{Discret}}$ to map $(\bm{x},\bm{P})$ into a chromosome.}
    \STATE{Adopt random scheme to obtain an initial population $\mathcal{U}$.}
     \STATE{Compute the fitness of each
chromosome in $\mathcal{U}$ by calling the solution in Section \ref{Subsec:TaskScheduingFixedPower}.}
		\WHILE {$l\le L$}
            \STATE {Replicate a new population $\mathcal{U}'$. }
		\STATE{Execute crossover and mutation operators on $\mathcal{U}'$.}
		\STATE{Map each chromosome in $\mathcal{U}'$ into $(\bm{x},\bm{P})$ and then compute their fitness values by calling the solution in Section \ref{Subsec:TaskScheduingFixedPower}. }
            \STATE{Merge the two populations of $\mathcal{U}$ and $\mathcal{U}'$.}
		\STATE{Use the strategies of tournament and elitism to select $U$ chromosomes from $\mathcal{U}\cup\mathcal{U}'$.}
		\STATE{$l = l+1$.}
		\ENDWHILE
		\STATE{Choose the best chromosome and then map it into $(\bm{x},\bm{P})$.}
		\ENSURE $(\bm{x},\bm{P})$.\\
	\end{algorithmic}
\end{algorithm}

In this section, we first use the result in Section \ref{Subsec:PowerAllocation} to reduce the solution space of  $\textbf{P0}(\bm{x},\bm{P})$ without losing optimality and then introduce our proposed two-layer solution.

We have shown that $\textbf{P1}(\bm{P})$ can be solved optimally in a closed form (i.e., $\bm{P}^*$). Therefore, we can obtain the minimum total energy consumption for ${\textbf{P0}(\bm{x},\bm{P}^*)}$. Furthermore, we substitute (\ref{Def:Rjk}) into (\ref{Def:pjk}) to obtain:
\begin{align}\label{Def:pjknew}
	p_{jk}=\frac{D_j}{B_c\log_2(1+\frac{P_{s(j)}G^{\text{tran}}_{s(j)}G^{\text{rec}}_hL_fL_l^{k}}{N})},~~\forall j\in\mathcal{J}, k\in\mathcal{K}_{s(j),\mathcal{H}},
\end{align}
which is inversely proportional to $P_{s(j)}$. This indicates that increasing the value $P_{s(j)}$ would shorten the transmitting time $p_{jk}$, further increasing the sum weights of tasks. As such, we can obtain that the optimal power allocation for $\textbf{P0}(\bm{x},\bm{P})$ is in the smaller set of feasible power, i.e., $\mathcal{P} = \{P_{s(j)}|P^*_{s(j)}\le P_{s(j)} \le P^{\text{max}}_{s(j)},~~\forall j\in \mathcal {J}\}$.

To tackle $\textbf{P0}(\bm{x},\bm{P})$ efficiently, we decompose it into the two levels of optimization as follows. In the upper-level optimization, we adopt a genetic framework to solve $\textbf{P0}(\bm{x},\bm{P})$ by optimizing variables $\bm{x}$ and $\bm{P}$. In the lower-level optimization, we use the proposed solution in Section \ref{Subsec:TaskScheduingFixedPower} to calculate the objective value of $\textbf{P0}(\bm{x},\bm{P})$ with given variables $\bm{x}$ and $\bm{P}$. To order to use the genetic framework to remodel $\textbf{P0}(\bm{x},\bm{P})$, we require the following two definitions.

\begin{figure}[t!]
    \centering
    \includegraphics[scale=0.45]{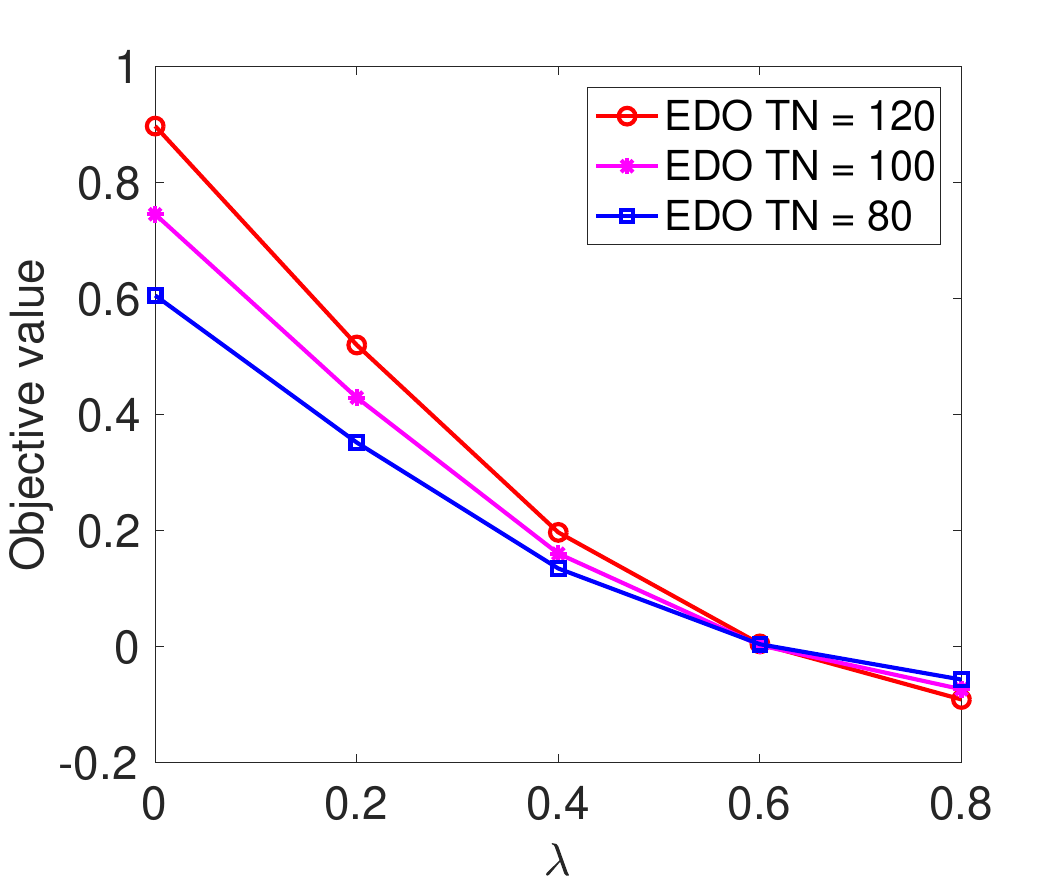}
    \caption{Objective value versus $\lambda$.}
    \label{fig2:ObjversusLambda}
\end{figure}
\begin{figure}[t!]
    \centering
    \includegraphics[scale=0.45]{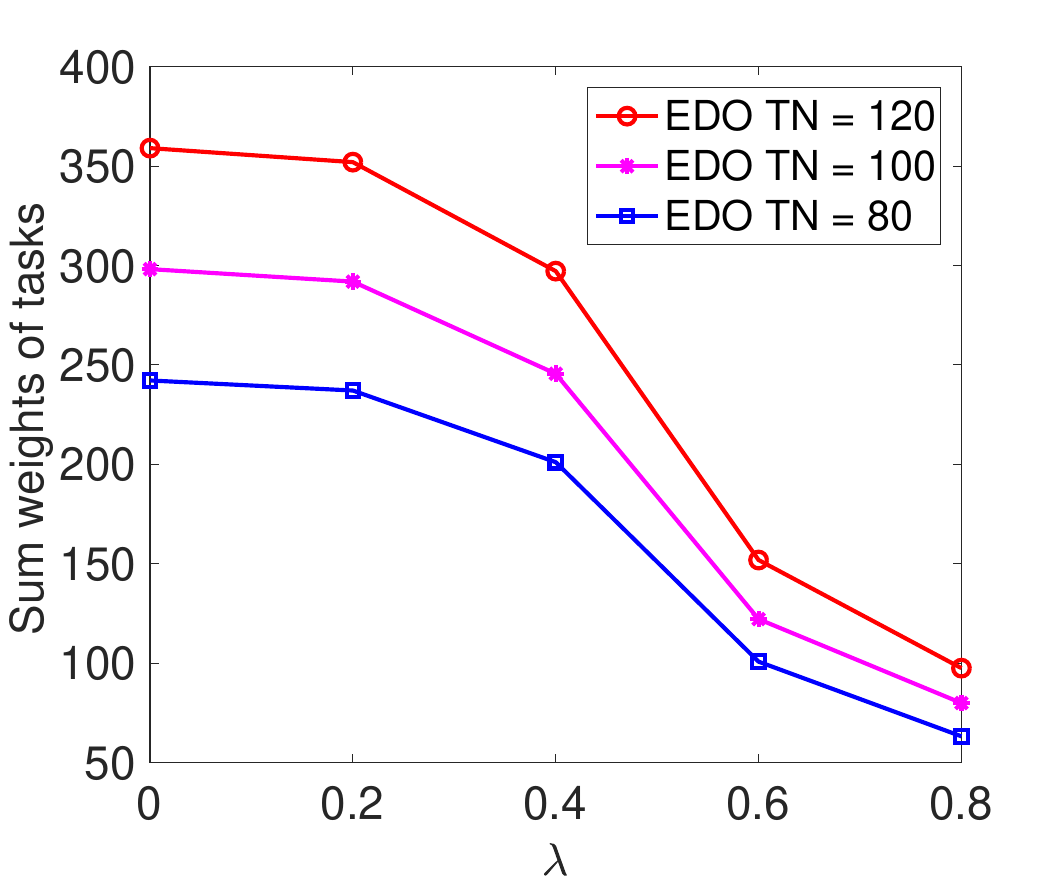}
    \caption{Sum weights of tasks versus $\lambda$.}
    \label{fig2:WVversusLambda}
\end{figure}
\subsubsection*{Genetic representation of the solutions to $\textbf{P0}(\bm{x},\bm{P})$} For convenient genetic representation, we first discretize the set of feasible power $\mathcal{P}$ to obtain a new set $\mathcal{P}_{\text{Discret}}=\{\mathcal{P}_{s(j)}\}$ with $\mathcal{P}_{s(j)} = \{P^k_{s(j)}| P^k_{s(j)}= P^*_{s(j)}+k\frac{P^{\text{max}}_{s(j)}-P^*_{s(j)}}{\left|\mathcal{P}_{s(j)}\right|}, k=0,1,...,\left|\mathcal{P}_{s(j)}\right|-1\}$. As such, we can represent each candidate discrete solution to $\textbf{P0}(\bm{x},\bm{P})$ as a chromosome. Specifically, we observe from the definition of $\bm{x}$ that each chromosome is actually a three-dimensional binary vector. Thus, we adopt binary coding to execute chromosome coding, such that one-to-one mapping relationship is established with $\left(\bm{x},\bm{P}\right)$ and a chromosome, by assigning zero or one to each element of each chromosome to indicate task scheduling and power allocation.
\subsubsection*{Quality evaluation of the represented solutions} We define a fitness function over the genetic representation as the objective value in $\textbf{P0}(\bm{x},\bm{P})$ to evaluate each candidate solution. The fitness of each chromosome can be computed as follows: First, we can obtain $\textbf{P4}(\bm{x})$ by substituting known $\bm{P}$ into $\textbf{P0}(\bm{x},\bm{P})$. Second, we can obtain the value of $\bm{x}$ from each chromosome through their one-to-one mapping relationship. Third, we use the value of $\bm{x}$ to construct a conflict graph $G(V,E)$ through the proposed method in Section \ref{Subsec:TaskScheduingFixedPower}. In particular, we represent each element of $\bm{x}$ equal to one as one vertex in $G(V,E)$ and compute a MWIS of $G(V,E)$. Then, we set the elements of $\bm{x}$ corresponding to the vertexes within MWIS to be ones; and zeros otherwise. Finally, we substitute both the obtained value of $\bm{x}$ and $\bm{P}$ into the objective function in $\textbf{P4}(\bm{x})$ to compute the fitness of each chromosome. Once the genetic representation and the fitness function are defined, Genetic Algorithm (GA) can begin with an initial population of chromosomes and then evolves better solutions to $\textbf{P0}(\bm{x},\bm{P})$ through the repetitive application of the bio-inspired operations of mutation, crossover, and selection.

\begin{figure}[t!]
    \centering
    \includegraphics[scale=0.45]{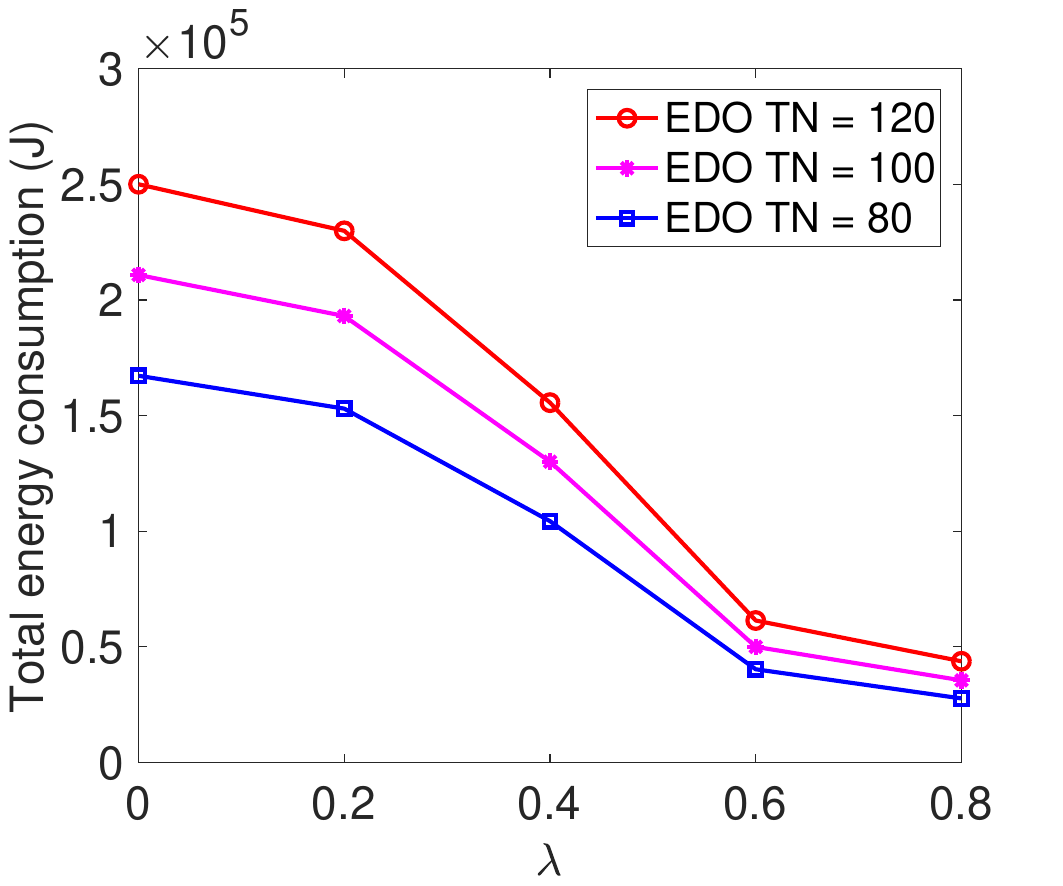}
    \caption{Total energy consumption versus $\lambda$.}
    \label{fig2:EEversusLambda}
\end{figure}
We term the above two-layer solution as EDO, which is summarized in Algorithm \ref{alg:EDO}. At the lower layer, solving $\textbf{P4}(\bm{x})$ has complexity $O(\left|V\right|^2)$ with $\left|V\right|$ indicating the number of vertexes in $G(V,E)$. At the higher layer, the complexity of GA is $O(UL)$, where $U$ and $L$ represent population size and generation number, respectively. Thus, the total complexity of EDO is $O(UL\left|V\right|^2)$.



\begin{figure}[t!]
    \centering
    \includegraphics[scale=0.45]{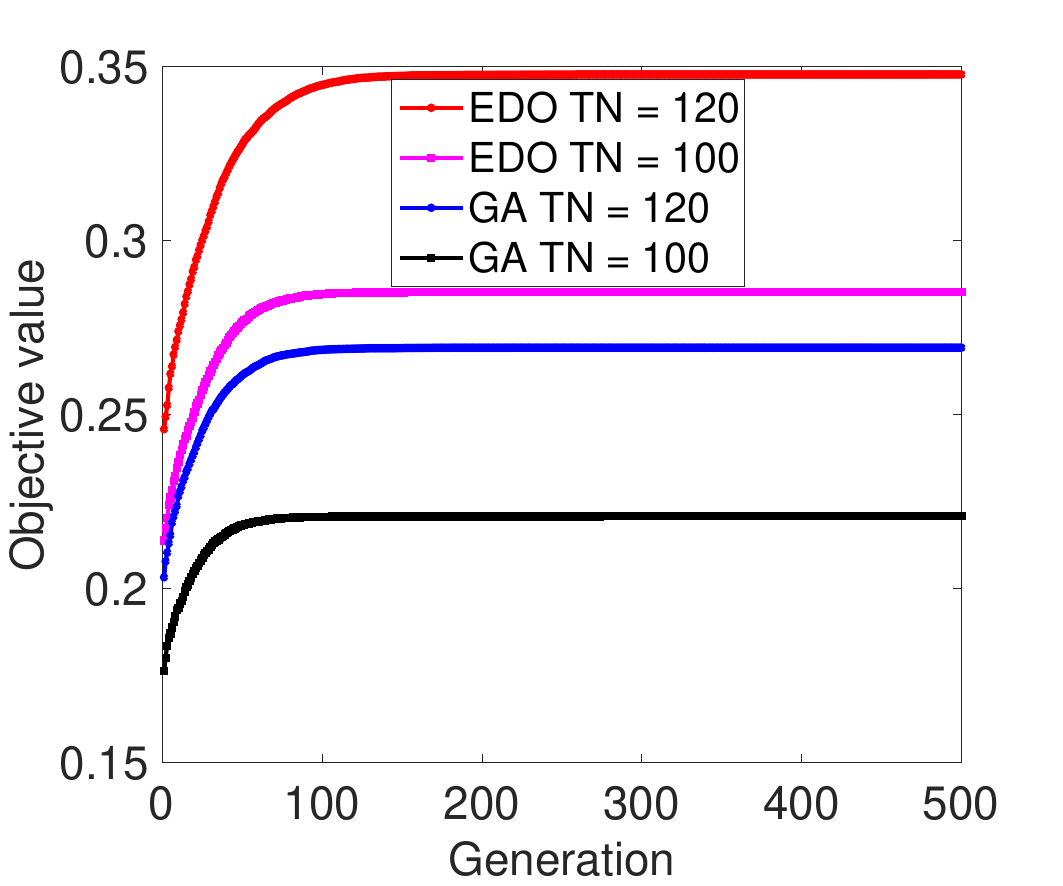}
    \caption{Objective value versus generations.}
    \label{fig1:Convergence}
\end{figure}

\section{Simulation results}
The simulation scenario consists of 4 GSs and 4 EOSs. We use the softwares of satellite tool kit (STK) and MATLAB to form a co-simulation platform to produce the considered simulation scenario. Specifically, by assigning the numbers of north American air (NORAD), we use MATLAB to input the two-line elements of each EOS into STK directly. The NORAD numbers of the four EOSs are set to 31113, 32382, 33320, and 32289. The locations of the four GSs are set to ($18^\circ\text{N}, 109.5^\circ \text{E}$), ($40^\circ\text{N}, 116^\circ \text{E}$), ($34^\circ\text{N}, 108^\circ \text{E}$), and ($39.5^\circ\text{N}, 76^\circ \text{E}$). Each GS or EOS is equipped with one antenna. The whole scheduling horizon is set to 12 hours. The size of time slot is set to 10 seconds. The data size of each task $j$ (i.e., $D_j$) is uniformly generated from the interval of $[500,1500]$ Mbits. The weight of each task $j$ (i.e., $w_j$) is randomly distributed within the interval of $[1,5]$. We set $G^{\text{tran}}_{s(j)}=G^{\text{rec}}_h= 36~\text{dB}$, $L_l^{k} = 1$, $L_f = 10^{-23}$, $N = 5.16*10^{-20}~\text{J}$, $B_c=2.2~\text{GHz}$, $R_{k}^{\text{req}} = 250~\text{Mbps}$, $P^{\text{max}}_{s(j)}=100~\text{Watts}$, $E_{\text{max}}^{\text{total}} = 2.5*10^5~\text{J}$, and $W=400$. The probabilities of crossover and mutation for GA are set to 0.5 and 0.8, respectively. We set $U=60$ and $L=200$.

\begin{figure}[t!]
    \centering
    \includegraphics[scale=0.45]{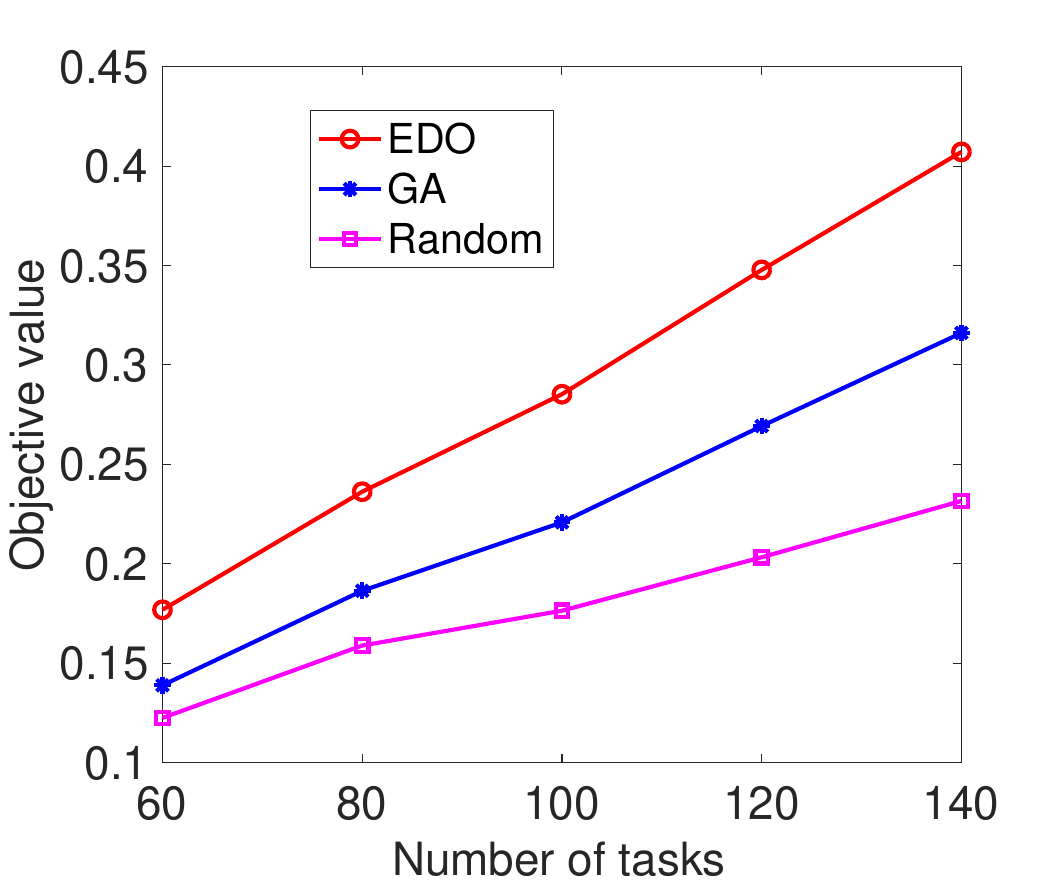}
    \caption{Objective value versus the number of tasks.}
    \label{fig3:ObjversusTN}
\end{figure}

We adopt two experiments to verify the performance of the proposed EDO algorithm. In the first experiment, we first show the design objective value versus $\lambda$ in Fig. \ref{fig2:ObjversusLambda}, and then vary the values of $\lambda$ to show the balance between two conflicting sub-objectives in Fig. \ref{fig2:WVversusLambda}, and Fig. \ref{fig2:EEversusLambda}. In the second experiment, we set $\lambda=0.3$ to respectively evaluate the convergence and superiority of EDO in Fig. \ref{fig1:Convergence} and Fig. \ref{fig3:ObjversusTN}. We adopt the baselines of \textit{GA} and \textit{Random} to evaluate the performance of EDO. In \textit{GA}, we use the solution in Section \ref{Subsec:TaskScheduingFixedPower} to calculate the fitness value for each chromosome and power allocation is subject to the interval of $[0,P^{\text{max}}_{s(j)}],\forall j$. In \textit{Random}, task scheduling and power allocation are randomly selected.

As depicted in Fig. \ref{fig2:ObjversusLambda}, we show how the objective value obtained by EDO varies with $\lambda$ for different task number (TN). From the objective value in Fig. \ref{fig2:ObjversusLambda}, we further plot its two terms of the sum weights of tasks and the total energy consumption in Fig. \ref{fig2:WVversusLambda} and Fig. \ref{fig2:EEversusLambda}, respectively. We can see that the objective value, the sum weights of tasks, and the total energy consumption decrease as $\lambda$ grows. This result coincides well with our design objective function in {\textbf{P0}}. This is because that a larger value of $\lambda$ favors the objective of minimizing the total energy consumption, while a smaller one favors the objective of maximizing the sum weights of tasks. This reveals a balance between the objectives of maximizing the sum weights of tasks and minimizing the total energy consumption. Furthermore, this phenomenon indicates that we can adjust the value of $\lambda$ in practical applications to produce the desired metrics involving the sum weights of tasks and the total energy consumption.

In Fig. \ref{fig1:Convergence}, we further show the comparisons between EDO and \textit{GA} in terms of the objective value versus generations. From Fig. \ref{fig1:Convergence}, we observe that the objective values obtained by both EDO and \textit{GA} first increase with the number of tasks and then remain constant after 150 generations. This demonstrates that the proposed EDO achieves fast convergence in limited generations and provides a better convergence performance over \textit{GA}.

In Fig. \ref{fig3:ObjversusTN}, we compare EDO with \textit{GA} and \textit{Random} in terms of the objective value versus the number of tasks. It is observed that the proposed EDO substantially outperforms the two comparison schemes in terms of the objective value as the number of tasks increase. This is because our proposed EDO jointly optimizes power allocation and task scheduling over a smaller feasible set of power to maximize the designed objective of {\textbf{P0}}. From Fig. \ref{fig3:ObjversusTN}, we can see that the objective value monotonically increases with the number of tasks increase. This is because that a smaller $\lambda$ tends to maximize the sum weights of tasks instead of minimizing the total energy consumption.

\section{Conclusion}
We study energy-efficient data offloading in EOSNs by jointly optimizing power allocation and task scheduling to minimize the total energy consumption and maximize the sum weights of tasks. We systematically explore the unique structures of joint optimization problem. We consider two special cases for the joint optimization problem and derive their near-optimal solutions. Then, we combine the two special solutions with genetic framework to propose an efficient two-layer EDO algorithm. Simulation results demonstrate the feasibility and effectiveness of the proposed algorithm.

\bibliographystyle{IEEEtran}
\bibliography{Reference}

\end{document}